\def\bbuildrel#1_#2^#3{\mathrel{\mathop{\kern 0pt#1}\limits_{#2}^{#3}}}

\def\NN{\mathbb N}

\def\QQ{\mathbb Q}
\def\PP{\mathbb P}

\def\RR{\mathbb R}

\def\Pot{\hbox{$\mathcal P$}}

\def\ms{\medskip}
\def\ss{\smallskip}

\def\w{\thinspace\hbox{\hsize 14pt \rightarrowfill}\thinspace}

\def\0{\hbox{$\emptyset$}}

\def\sub{\subseteq}

\def\K{\mathscr{K}}

\documentclass[a4paper,12pt]{amsart}

\usepackage{amssymb,enumerate,mathrsfs, color}
\usepackage[utf8]{inputenc} 

\theoremstyle{plain}

\newcommand{\co}{\mathfrak c}

\newcommand{\cantor}{2^{\NN}}
\newcommand{\baire}{\NN^{\NN}}
\newtheorem{theorem}{Theorem}[section]

\newtheorem{lemma}[theorem]{Lemma}

\newtheorem{definition}[theorem]{Definition}

\newtheorem{proposition}[theorem]{Proposition}

\newtheorem{remark}[theorem]{Remark}
\newtheorem{cor}[theorem]{Corollary}

\numberwithin{equation}{section}

\begin{document}

\title{A note on uniform continuity of monotone 
	 functions}

\author{Roman Pol, Piotr Zakrzewski and Lyubomyr Zdomskyy}


\address{Roman Pol, Institute of Mathematics, University of Warsaw,  Banacha 2, 02-097 Warsaw, Poland}
\email{pol@mimuw.edu.pl}

\address{Piotr Zakrzewski, Institute of Mathematics, University of Warsaw,  Banacha 2, 02-097 Warsaw, Poland}
\email{piotrzak@mimuw.edu.pl}

\address{Lyubomyr Zdomskyy, Institut f\" ur Diskrete Mathematik und Geometrie, Technische Universit\" at Wien,  Wiedner Hauptstrasse 8-10/104, 1040 Wien, Austria}
\email{lzdomskyy@gmail.com}

\makeatletter
\@namedef{subjclassname@2020}{\textup{2020} Mathematics Subject Classification}
\makeatother



\thanks{The  research of the third named author was funded in whole by the Austrian Science Fund (FWF), Grants DOI 10.55776/I5930 and 
10.55776/PAT5730424.}


\begin{abstract}
	 We prove that it is consistent with ZFC that
	for every non-decreasing
	function $f:[0,1]\w [0,1]$, each subset of $[0,1]$ of cardinality $\co$ contains a set of cardinality $\co$ on which $f$ is uniformly continuous. We show that this statement follows from the assumptions that $\mathfrak d^* <\co$ and $\co$ is regular, where $\mathfrak d^*\leq \mathfrak d$ is the smallest cardinality $\kappa$ such that any two disjoint countable  dense sets  in  $\cantor$  can be separated by sets each of which is an intersection of at most $\kappa$-many open sets in $\cantor$. 
	 We establish also that  $\mathfrak d^*=\min\{\mathfrak u, \mathfrak d\}=\min\{\mathfrak r, \mathfrak d\}$, thus giving an alternative proof of the latter equality established by J. Aubrey in 2004.

\end{abstract}

\maketitle

\section{Introduction}

The main {subject} of this note is the following theorem.  

\begin{theorem}\label{con}
	The statement
	\begin{itemize}
		\item[(*)] For every non-decreasing
		function $f:[0,1]\w [0,1]$, each subset of $[0,1]$ of cardinality $\co$ contains a set of cardinality $\co$ on which $f$ is uniformly continuous.
	\end{itemize}
	is independent of ZFC.
	
\end{theorem} 

Sierpiński \cite[Th\'{e}or\`{e}me 6]{Si} proved that CH implies the negation of $(*)$. Sierpiński considered 
an increasing function $f:[0,1]\w [0,1]$ which is discontinuous precisely  at the rationals in $(0,1)$ (a well-known example of such a function is due to Lebesgue and is defined by letting $L(0)=0$ and $L(x)=\sum_{\{n\in\NN: \ q_n<x\}} 2^{-n}$ for $x\in (0,1]$, where $\{q_n:\ n\in \NN\}$ is an injective enumeration of $\QQ\cap [0,1)$). He
proved that the (continuous) restriction of $f$ to the set $\PP$ of  irrationals in $[0,1]$ is not uniformly continuous on any Lusin set in $\PP$.

The fact that in some models of ZFC statement $(*)$ can be true will be explained in Section \ref{Section 2}, where a new cardinal coefficient $\mathfrak d^*$ which plays a crucial role in our considerations is introduced. A combinatorial analysis, giving a precise description of $\mathfrak d^*$ in terms of some other known cardinal characteristics of the continuum, will be presented in Section \ref{Section 3}. Some additional comments will be gathered in Section \ref{Section 4}.

\section{The statement $(*)$ is consistently true}\label{Section 2}


\begin{definition}

		Let $\mathsf{Fin}$ and $\mathsf{cFin}$ be the families of all finite and co-finite subsets of $\mathbb N$, respectively. Then
		$\mathfrak d^*$ is the minimal cardinality of a cover $\mathsf K$
		of $\mathcal P(\mathbb N)$  by compact subspaces\footnote{As usual,
			we identify $\mathcal P(\mathbb N)$ with $2^{\mathbb N}$ using characteristic functions.} such that for each $\mathcal K\in\mathsf K$,
		either $\mathcal K\cap\mathsf{Fin}=\emptyset$ or $\mathcal K\cap\mathsf{cFin}=\emptyset$.
	
\end{definition}

Since both $\mathcal P(\mathbb N)\setminus\mathsf{Fin}$ and
$\mathcal P(\mathbb N)\setminus\mathsf{cFin}$ are homeomorphic to the Baire space
$\mathbb N^{\mathbb N}$ which can be covered by $\mathfrak d$ many compact subspaces, we have $\mathfrak d^*\leq \mathfrak d$ (here $\mathfrak{d}$ denotes, as usual, the {\sl dominating number}, i.e., the smallest cardinality of a dominating family in $\baire$ corresponding to the ordering of eventual domination $\leq^*$, cf. \cite{Bla10}). It follows that it is consistent with ZFC to assume that $\mathfrak d^* <\co$ and $\co$ is regular, cf. \cite{Bla10}, and we shall show that these conditions yield $(*)$.

\ms 

Let us start from the following observation.

\begin{lemma}\label{compact}
	Given two disjoint countable sets $D_0, D_1$ in a compact metrizable space $X$, there is a collection $\K$ of compact sets in $X$ such that $|\K|\leq \mathfrak d^*$, $\bigcup\K=X$ and each element of $\K$ hits at most one $D_i$.
\end{lemma}

\begin{proof}
	 Using the Cantor–Bendixson theorem we write $X$
	in the form $Y\cup Z$, $Y\cap Z=\emptyset$, such that $Y$ is a compact space without isolated points and $Z$ is countable. Now, to get $\K$,  it suffices to find a relevant cover of $Y$  and then to extend it by the singletons of $Z$. Instead, we simply assume that $X$ has no isolated points.
	
	\ss 
	
	 There is a continuous surjection
	 $u:\Pot(\NN) \to X$
	 such that $u^{-1}(d)$ is a singleton, for $d\in D_0\cup D_1$. Indeed, if we expand each $d\in D=D_0\cup D_1$ to a copy $P_d$ of the Cantor set in $X$, then \cite[\S 45, II, Theorem 4]{Ku} gives us   a continuous surjection
	 $u:\Pot(\NN) \to X$ which is injective on $u^{-1}(\bigcup_{d\in D}P_d)$ . 
	
	\ss 
	
	Now, the disjoint sets $u^{-1}(D_i)$, $i=0,1$, are countable, and by adopting standard arguments concerning countable dense homogeneity of $\cantor$ (cf. \cite[Theorem 1.6.9]{vM}) one can find a homeomorphism 
	$h:\Pot(\NN)\to\Pot(\NN)$
	 such that $u^{-1}(D_0)\sub h(\mathsf{Fin})$ and $u^{-1}(D_1)\sub h(\mathsf{cFin})$. 
	 
	 \ss 
	 
	 Finally, if $\mathsf K$ is a collection of compact sets in $\Pot(\NN)$ satisfying the conditions in the definition of $\mathfrak d^*$, with $|\mathsf K|=\mathfrak d^*$, the collection $\K=\{(u\circ h)({\mathcal K}):\ {\mathcal K}\in\mathsf K\}$ satisfies the assertion of the lemma. 
	
\end{proof}

To complete  the proof of the consistency of statement $(*)$ (hence also the proof of Theorem \ref{con}) it is enough to prove the following result.

\begin{proposition}\label{main}
	If $\mathfrak d^* <\co$ and $\co$ is regular, then statement $(*)$ is true.
\end{proposition} 

\begin{proof} Let
	a function $f:[0,1]\w [0,1]$ be	non-decreasing. 
	Let $D$ be the set of discontinuity points of $f$ and let X be the closure in the square $[0,1]\times [0,1]$ of the graph of $f$ restricted to the set
	$[0,1] \setminus D$. 
	 For each $t\in [0,1]$ let $X(t)=X\cap (\{t\}\times [0,1])$.
	 
	  Since $f$ is non-decreasing, the set $D$ of discontinuity points of $f$ is countable and  for each $d\in D$ we have $X(d)=\{(d,d_0),(d,d_1)\}$, where $d_0< d_1$, while $X(t)=\{(t,f(t))\}$ for $t\in [0,1]\setminus D$.
	
	Let 
	$$
	Q_i=\{(d,d_i):\ d\in D\},\ i=0,1.
	$$
	
	Since $X$ is compact, Lemma \ref{compact} provides a covering  $\K$ of $X$ by  compact sets such that  $|\K|\leq\mathfrak d^*$ and each element
	of $\K$ intersects at most one of the sets  $Q_0,\ Q_1$. 
	
	Let $\pi: X\w [0,1]$ be the projection onto the first coordinate. Since each $K\in\K$ contains at most one point from every pair $X(d)=\{(d,d_0),(d,d_1)\}$, $d\in D\cap \pi(K)$, $K$ is the graph of a continuous function $f_K: \pi(K)\w [0,1]$. Since $K$ is compact, $f_K$ is uniformly continuous. Moreover, $f$ and $f_K$ coincide on $\pi(K)\setminus D$,  so $f$ is uniformly continuous on $\pi(K)\setminus D$.
	
	We are ready to address statement $(*)$. Let $E\sub [0,1]$ have cardinality $\co$. Since the sets $\pi(K)$, $K\in\K$, cover $[0,1]$ and $|\K|\leq\mathfrak d^* <\co$ (cf. Lemma \ref{compact}), the regularity of $\co$ implies that there is $K\in\K$ with $|\pi(K)\cap E|=\co$, and since $f$ is uniformly continuous on the set  $(\pi(K)\cap E)\setminus D$, we get $(*)$. 
\end{proof}

\section{$\mathfrak d^*=\min\{\mathfrak u, \mathfrak d\}=\min\{\mathfrak r, \mathfrak d\}$}\label{Section 3}

The aim of this section it to prove the equalities announced in its title.
Let us recall that the {\sl ultrafilter number} $\mathfrak{u}$ is the minimal size of a base of a free ultrafilter in $\Pot(\NN)$, and the {\sl reaping number}
$\mathfrak r$ is the minimal size of a {\sl reaping subfamily} $\mathcal B\sub\Pot(\NN)$, i.e., a family $\mathcal B$ of infinite subsets of $\NN$ such that for every $X\sub\NN$ there exists $B\in\mathcal B$
such that either $X\supset B$ or $B\cap X=\emptyset$,
see \cite{Bla10} for more information on these cardinals.

Clearly,
each ultrafilter base is  reaping, hence $\mathfrak r\leq\mathfrak u$,
and the strict inequality is consistent, cf. \cite{GolShe90}. However,
by \cite[Corollary 6.3]{Aub04} we cannot have
 $\mathfrak r< \mathfrak u$  if $\mathfrak r<\mathfrak d$, which is clearly equivalent to $\min\{\mathfrak u, \mathfrak d\}=\min\{\mathfrak r, \mathfrak d\}$. In the course of our proof 
 of Theorem~\ref{d^* versus d,u} below we reestablish the latter 
 equality, this way giving an alternative and more streamlined proof thereof, as well as giving a natural topological interpretation of
 $\min\{\mathfrak u, \mathfrak d\}$.

\begin{definition}\hfill\null
\begin{itemize}
\item A family $\mathcal F$ of infinite subsets of $\mathbb N$
is called a \emph{semifilter}, if it is closed under taking {almost-supersets} of its elements, i.e., if
$F\in \mathcal F$ and $F\subseteq^*X\subseteq\mathbb N$, then $X\in\mathcal F$.
\item A family $\mathcal C$ of infinite subsets of $\NN$ is called
\emph{centered} if $\cap\mathcal C'$ is infinite for any finite subfamily $\mathcal C'\subseteq\mathcal C$.
\end{itemize}
\end{definition}

For a semifilter $\mathcal F$ we denote by $\mathcal F^+$ the family 
$\{X\subseteq\mathbb N:\forall F\in\mathcal F\:(F\cap X\neq\emptyset)\}$.
It is clear that $\mathcal F^+$ is also a semifilter.

For $\mathcal X \subseteq\mathcal P(\mathbb N)$ we denote by
$\sim\!\mathcal X$ the family $\{\mathbb N\setminus X:X\in\mathcal X\}$.
It is easy to check that $\mathcal F^+=\mathcal P(\mathbb N)
\setminus\sim\!\mathcal F$. Note that if $\mathcal F_0\subseteq\mathcal F_1$ are semifilters,
then $\mathcal F_0^+\supseteq\mathcal F_1^+$.

\ss 

The following two facts are probably well-known, but we give their proofs 
for the sake of completeness. The first one is remiscent of 
\cite[Theorem~2.10]{Bla10}.

\begin{lemma}\label{a_la_blass}
Let $\kappa<\mathfrak d$ be a cardinal and
$$\big\{ \langle n^{\alpha}_{i}:i\in\NN\rangle\: :\: \alpha<\kappa\big\}$$
be a family of increasing number sequences. Then there exists an increasing number sequence $\langle h(j):j\in\NN\rangle$
such that for every finite $A\subset\kappa$ the set
$$J_A=\{j\in\NN\: :\: \forall \alpha\in A\, \exists i\in\NN\:\big([n^{\alpha}_{i}, n^{\alpha}_{i+1})\sub [h(j), h(j+1) )\big) \}$$
is infinite.
\end{lemma}
\begin{proof}
It is enough to find $h$ satisfying the statement of the lemma for
singletons $A\sub\kappa$, i.e., $A$ of the form $\{\alpha\}$ for some 
$\alpha\in\kappa$, since without loss of generality we may assume that for every finite $A\sub\kappa$ there exists $\gamma\in\kappa$
such that 
$$ \forall j\in\NN\, \forall \alpha\in A\,\exists i\in\NN\:
\big([n^{\alpha}_{i}, n^{\alpha}_{i+1})\sub[n^{\gamma}_{j}, n^{\gamma}_{j+1}) \big).$$
For every $\alpha\in\kappa$ and $i\in\NN$ set $f_\alpha(i)=n^\alpha_{2i}$
and find a strictly increasing $h\in\NN^\NN$ such that $h(0)=0$ and
$h\not\leq^* f_\alpha$ for all $\alpha\in\kappa$. We claim that 
$\langle h(j):j\in\NN\rangle$ is as required. Indeed, let us fix  $i>0$ such that 
$h(i)\geq f_\alpha(i)=n^\alpha_{2i}$. Thus,
$$ \{n^\alpha_k:k\leq 2i\}\sub \bigcup_{l<i}[h(l),h(l+1)), $$
which gives that 
$ \big|\big\{k\ { < }\ 2i : n^\alpha_k,n^\alpha_{k+1}$ belong to the same interval
$ [h(l),h(l+1)) \big\}\big|\geq i $, and consequently 
$J_{\{\alpha\}}$ defined in the formulation must be infinite because 
$i$ with $h(i)\geq f_\alpha(i)=n^\alpha_{2i}$ may be taken arbitrarily large.
\end{proof}

\begin{cor}\label{less_d}
Let $\mathsf K$ be a family of $<\mathfrak d$-many compact subspaces
of $\mathcal P(\mathbb N)\setminus\mathsf{Fin} $.
Then there exists a monotone surjection $\phi:\mathbb N\to\mathbb N$
and a centered family $\mathcal C$ of size $|\mathcal C|\leq |\mathsf K|$
such that 
$\{\phi[X]: X \in \bigcup K\}\subseteq\mathcal F$ , where 
$\mathcal F$ is the smallest free filter generated by $\mathcal C$.
%
%
%
\end{cor}
\begin{proof}
For every $\mathcal K\in\mathsf K$, $\mathcal K$ being compact, there is an increasing number sequence $\langle n^{\mathcal K}_{i}:i\in\NN\rangle$  
such that $X\cap [n^{\mathcal K}_{i},n^{\mathcal K}_{i+1})\neq\emptyset$
for any $X\in\mathcal K$ and $i\in\mathbb N$.
Since $|\mathsf K|<\mathfrak d$, by 
Lemma~\ref{a_la_blass}
there exists a strictly increasing function 
$h:\NN\to\NN$ such that $h(0)=0$ and for any finite $\mathsf K'\subseteq\mathsf K$
the set $C(\mathsf K')$ of those $j\in\NN$ such that for every $\mathcal K\in\mathsf K'$ there exists $i\in\NN$ with 
\begin{equation}\label{gr_abstand}
[h(j),h(j+1))\supseteq [n^{\mathcal K}_{i},n^{\mathcal K}_{i+1}),
\end{equation}
is infinite. Note that  the family
 $$\mathcal C=\{C(\mathsf K'):\mathsf K'\mbox{ is a finite subset of }\mathsf K\}$$
is centered because by the definition of $C(\mathsf K')$ we have
$$C(\mathsf K'_0\cup \mathsf K'_1\cup\cdots\cup \mathsf K'_n)=\bigcap_{i\leq n}C(\mathsf K'_i)$$
for any finite family $\{\mathsf K'_0,\mathsf K'_1,\ldots,\mathsf K'_n\}$
of finite subsets of $\mathsf K$.

We claim that $\mathcal C$ and $\phi:\NN\to\NN$  such that $\phi^{-1}(j)=[h(j),h(j+1))$ for all $j\in\NN$, are as required. Indeed, given 
$X\in\bigcup\mathsf K$, let $\mathcal K\in\mathsf K$ be such that $X\in\mathcal K$.
Then 
\begin{eqnarray*}
\phi[X]=\{j\in\NN: X\cap [h(j),h(j+1))\neq\emptyset\}\supseteq \\
\supseteq
\big\{j\in\NN\: :\: \exists i\in\NN\: \big([h(j),h(j+1))\supseteq [n^{\mathcal K}_{i},n^{\mathcal K}_{i+1})\big)\big\}\supseteq C(\{\mathcal K\})\in\mathcal C,
\end{eqnarray*}
which completes the proof.
\end{proof}

For every $X\subseteq\mathbb N$ we adopt the following notation:
\begin{eqnarray*}
\downarrow\!X =\mathcal P(X), \ \ 
 \uparrow\!X = \{Y\subseteq\mathbb N:X\sub Y\},\\
 \downarrow^*\!X = \{Y\subseteq\mathbb N:Y\subseteq^* X\},\ \ 
 \uparrow^*\!X = \{Y\subseteq\mathbb N:X\subseteq^* Y\}.
\end{eqnarray*}

 For 
a subset $\mathcal L\subseteq\mathcal P(\mathbb N)$ we set: 
\begin{eqnarray*}
\uparrow\!\mathcal L=\bigcup\{\uparrow\!X:X\in \mathcal L\},\ \ 
\downarrow\!\mathcal L=\bigcup\{\downarrow\!X:X\in \mathcal L\},\\
\uparrow^*\!\mathcal L=\bigcup\{\uparrow^*\!X:X\in \mathcal L\},\ \ 
\downarrow^*\!\mathcal L=\bigcup\{\downarrow^*\!X:X\in \mathcal L\}. 
\end{eqnarray*}

 It is easy to see that if $\mathcal L$ is compact,
then both $\uparrow\!\mathcal L$ and $\downarrow\!\mathcal L$ are compact, and 
$\uparrow^*\!\mathcal L$ and $\downarrow^*\!\mathcal L$ are $\sigma$-compact.

\ss 

We are now ready to prove the main result of this section.

\begin{theorem}\label{d^* versus d,u}
$\mathfrak d^*=\min\{\mathfrak u,\mathfrak d\}=\min\{\mathfrak r,\mathfrak d\}$.
\end{theorem}
\begin{proof}
As we have already noticed, since both $\mathcal P(\mathbb N)\setminus\mathsf{Fin}$ and
$\mathcal P(\mathbb N)\setminus\mathsf{cFin}$ are homeomorphic to the Baire space
$\mathbb N^{\mathbb N}$ which can be covered by $\mathfrak d$ many compact subspaces, we have $\mathfrak d^*\leq \mathfrak d$.

Let $\mathcal B$ be a reaping family
of cardinality
$\mathfrak r$. For every $B\in\mathcal B$ set  $\mathcal K(B)_d=\downarrow\!(\mathbb N\setminus B)$ and $\mathcal K(B)_u=\uparrow\!B$.

It follows that 
$$\mathcal P(\mathbb N)=\bigcup\{\mathcal K(B)_d:B\in\mathcal B\}\cup\bigcup \{\mathcal K(B)_u:B\in\mathcal B\}$$
and 
$\mathcal K(B)_d\cap \mathsf{cFin}=\mathcal K(B)_u\cap \mathsf{Fin}=\emptyset$ for all $B\in\mathcal B$,
 hence $\mathfrak d^*\leq\mathfrak r$, which together with the first paragraph yields $\mathfrak d^*\leq \min\{\mathfrak r,\mathfrak d\}$.
\smallskip

To show that $\mathfrak d^*\geq \min\{\mathfrak u,\mathfrak d\}$, suppose, towards a contradiction, that there is a family $\mathsf K=\mathsf K_d\cup\mathsf K_u$ 
of compact subspaces of $\mathcal P(\mathbb N)$ of size  $\kappa=|\mathsf K|<\min\{\mathfrak u,\mathfrak d\}$ such that $\bigcup\mathsf K=\mathcal P(\NN)$ and
\begin{equation}\label{eq_01}
\mathcal K_0\cap\mathsf{cFin}=\mathcal K_1\cap\mathsf{Fin}=\emptyset
\end{equation}
for 
any $\mathcal K_0\in \mathsf K_d$ and $\mathcal K_1\in \mathsf K_u$.

Replacing each $\mathcal K\in\mathsf K_d$  (resp. $\mathcal K\in\mathsf K_u$)
with a countable family of compact subspaces of $\mathcal P(\NN)$ covering $\downarrow^*\!\mathcal K$ (resp. $\uparrow^*\!\mathcal K$),
we can assume that $\mathcal U_u=\bigcup\mathsf K_u$ as well as 
$\mathcal U_d=\sim\!\bigcup\mathsf K_d$ are  semifilters on $\NN$.
Set $\mathcal U=\mathcal U_u\cup\mathcal U_d$.
It follows that 
$$\mathcal U_d^+=\mathcal P(\mathbb N)\setminus\sim\!\mathcal U_d=
\mathcal P(\mathbb N)\setminus\sim\!(\sim\!\bigcup\mathsf K_d)=
\mathcal P(\mathbb N)\setminus \bigcup\mathsf K_d\sub \bigcup\mathsf K_u=\mathcal U_u,$$
and hence 
\begin{equation} \label{big_semif}
\mathcal U^+\sub \mathcal U_d^+\sub \mathcal U_u\sub \mathcal U.
\end{equation}

By the construction, the semifilter $\mathcal U$ can be covered by $\kappa$
many compact subspaces. Since $\kappa<\mathfrak d$,
by Corollary~\ref{less_d}, there exists a free filter $\mathcal F$ generated 
by at most $\kappa$ many subsets of $\NN$ and a monotone surjection $\phi:\NN\to\NN$
with $\mathcal F\supseteq
\{\phi[X]: X \in \mathcal U\}$. Moreover, the family $\{\phi[X]: X \in \mathcal U\}$ is a semifilter and one easily checks that
\begin{eqnarray*}
\mathcal F^+\sub \{\phi[X]: X \in \mathcal U\}^+ =
\{\phi[X]: X \in \mathcal U^+\}\sub \{\phi[X]: X \in \mathcal U\}\subseteq\mathcal F,
\end{eqnarray*}
the second inclusion being a consequence of (\ref{big_semif}). This shows that
$\mathcal F$ is an ultrafilter generated by at most $\kappa$ many sets,
which is impossible because $\kappa<\mathfrak u$.
\smallskip

It follows from the above that 
$$\min\{\mathfrak u,\mathfrak d\}\leq \mathfrak d^* \leq \min\{\mathfrak r,\mathfrak d\},$$
and therefore these two inequalities must actually be equalities 
because $\mathfrak r\leq\mathfrak u$. 
\end{proof}

\section{Comments}\label{Section 4}

\subsection{Non-decreasing functions versus arbitrary functions   }\label{sec:4.1}

In \cite[Poposition 3.5]{P-Z} the following result is obtained:

\ms 
{\sl Let ${g}:\PP\w [0,1]$ be a continuous function such that the closure of the graph of $f$ hits each section $\{q\}\times [0,1]$, $q\in \QQ\cap [0,1]$, in an uncountable set. Then $\PP$ cannot be covered by less than $\mathfrak{d}$ sets on which ${g}$ is uniformly continuous. In particular, if $\mathfrak{d}=\co$, there is a subset $E$ of $\PP$ of cardinality $\co$ such that ${g}$ is not uniformly continuous on any subset of $E$ of cardinality $\co$} (cf. the proof of \cite[Theorem 3.4]{P-Z}).

\ms 
{In contrast, Proposition \ref{main} yields the following.}

\begin{remark}\label{generalization}
 {Let $h:E \to [0,1]$ be a monotone function  
	 	 on a set $E\sub [0,1]$ of cardinality $\co$. 
	 	  If  $\mathfrak d^* <\co$ and $\co$ is regular, then $h$ is uniformly continuous on
	  	  a set of cardinality $\co$.}

  \end{remark}
  
  {Indeed, if $h:E\to[0,1]$ is, say, non-decreasing,
  then we can extend it to $f:[0,1]\to [0,1]$
  which is also non-decreasing (by setting $f(y)=\sup h[E\cap [0,y] ]$
  for all $y\in [0,1]$, where we set $\sup(\emptyset)=0$)
  and then we can apply Proposition \ref{main}.}
	  
	  \ss 
	  
	 {Remark \ref{generalization}} applied to the continuous function ${g}$ and the set $E$ described in the opening statement of this subsection, leads to the following.

	  {
	  		\begin{remark}
	  			If $\mathfrak d^* <\mathfrak d=\co$ and $\co$ is regular, then there exists a continuous function $f: E\to [0,1]$ on a set $E\sub\RR$ of cardinality $\co$ such that no restriction of $f$ to a subset of $E$  of cardinality $\co$ is monotone.
	  	\end{remark}
  	}


\subsection{Models of ZFC with $\mathfrak d^* <\mathfrak d$}\label{sec:4.2} 

In view of the remarks from the preceding subsection and for the sake of completeness, we would like to mention that an example of a model of ZFC with $\mathfrak d^*=\mathfrak u = \aleph_1 < \aleph_2 = \mathfrak d$ (cf. Theorem \ref{d^* versus d,u})  is provided by the Miller model resulting from the $\aleph_2$-length countable support iteration of Miller forcing over a model of GCH, see \cite[Model 7.5.2]{Ba-Ju}. Another model which has the advantage
that $\mathfrak d^*=\mathfrak u < \mathfrak d$ can be any prescribed uncountable regular cardinals is presented by Blass and Shelah in \cite{Bl-Sh}.

\subsection{A generalization of $\mathfrak d^*$}
In view of Lemma \ref{compact} it is natural to generalize $\mathfrak d^*$ as follows.
	\begin{definition}
	For a  compact  metrizable space $X$ and countable disjoint sets
	$D_0,\ D_1$ in $X$ let 
	$\mathfrak d^*(X,D_0,D_1)$ be the minimal cardinality of a cover $\mathsf K$
	of $X$  by compact subspaces such that for each $\mathcal K\in\mathsf K$,
	either $\mathcal K\cap D_0=\emptyset$ or $\mathcal K\cap D_1=\emptyset$.
\end{definition}

In this notation we have $\mathfrak d^*=\mathfrak d^*(\Pot(\NN),\mathsf{Fin}, \mathsf{cFin})$ and Lemma  \ref{compact} states that
$\mathfrak d^*(X,D_0,D_1)\leq \mathfrak d^*$ for any
 two disjoint countable sets $D_0, D_1$ in an arbitrary compact metrizable space $X$. We complement this lemma with the following observation (whose proof gives also another method for establishing   Lemma  \ref{compact}).

\begin{proposition}\label{gen}
	Let $X$ be	a  compact  metrizable space with no isolated points and let $D_0,\ D_1$  be two disjoint countable  dense subspaces of $X$. Then 
	$\mathfrak d^*(X,D_0,D_1)=\mathfrak d^*$ (equivalently, $\mathfrak d^*$ is the smallest cardinality $\kappa$ such that $D_0,\ D_1$ can be separated by two sets each of which is the intersection of a collection of at most $\kappa$-many open sets in $X$).
\end{proposition}

In the proof we shall need the following well-known fact which can be proved by the standard back-and-forth argument.

\begin{lemma}\label{lem_q}
	Let $\mathbb Q$ be the set of the rational numbers and 
	$\mathbb Q=Q_0^0\cup Q_1^0$, $\mathbb Q=Q_0^1\cup Q_1^1$
	be two decompositions of $\mathbb Q$ into disjoint dense subspaces.
	Then there exists an automorphism of $\langle\QQ,\leq\rangle$ (hence a homeomorphism) $h:\mathbb Q\to\mathbb Q$
	such that $h[Q_0^0]=Q^1_0$ and $h[Q_1^0]=Q^1_1$.
\end{lemma}

\begin{proof}[Proof of Proposition \ref{gen}]
As	$X$ has no isolated points and $D_0,\ D_1$ are dense in $X$,
	 both $A=\mathsf{Fin}\cup\mathsf{cFin}\subseteq\mathcal
	P(\NN)$
	and $B=D_0\cup D_1$ are homeomorphic to the rationals
	$\mathbb Q$ being countable metrizable spaces without isolated points. Let $h:A\to B$ be a homeomorphism
	such that $h[\mathsf{Fin}]=D_0$ and $h[\mathsf{cFin}]=D_1$
	(its existence follows from Lemma~\ref{lem_q})
	and let $\beta h: \beta A\to\beta B$ be its homeomorphic extension to the \v{C}ech-Stone compactification $ \beta A$ of $A$.
	
	By \cite[Theorem~3.5.7]{eng-2} we have $\beta h[\beta A\setminus A]=\beta B\setminus B$ and analogously, the  extensions $f_A:\beta A\to \mathcal P(\NN)$ and $f_B:\beta B\to X$ of the identity maps $i_A:A\to \Pot(\NN)$ and
	$i_B:B\to X$ have the  property that $f_A[\beta A\setminus A]=\mathcal P(\NN)\setminus A$ and $f_B[\beta B\setminus B]=X\setminus B$.
	
	 It follows that if $\mathsf K$ is a family of compact subspaces of $\mathcal P(\NN)$ such as in the definition 
	of $\mathfrak d^*$, then 
	$$\{f_B[\beta h[f_A^{-1}[\mathcal K]]]:\mathcal K\in\mathsf K\}$$
	is such as in the definition 
	of $\mathfrak d^*(X, D_0,D_1)$; and vice versa,
	if $\mathsf K$ is a family of compact subspaces of $X$ such as in the definition 
	of $\mathfrak d^*(X, D_0,D_1)$, then 
	$$\{f_A[\beta h^{-1}[f_B^{-1}[\mathcal K]]]:\mathcal K\in\mathsf K\}$$
	is such as in the definition 
	of $\mathfrak d^*$. Thus, $\mathfrak d^*(X, D_0,D_1)=\mathfrak  d^*$, which completes the proof.
	
\end{proof}

\bibliographystyle{amsplain}

\end{document}